\documentstyle[12pt]{article}
\font\erl=eufm10 at 12pt
\newcommand{\euf}[3]{\mbox{\erl \char'#1#2#3}}
\font\ers=eufm10 

\font\erb=eufm10 at 15pt
\newcommand{\eufb}[3]{\mbox{\erb \char'#1#2#3}}
\font\ers=eusm10 at 12pt
\newcommand{\eus}[3]{\mbox{\ers \char'#1#2#3}}

\font\ms=msam8 at 11pt
\newcommand{\eql}{\mbox{\thinspace\thinspace{\ms 5}\thinspace\thinspace}}
\newcommand{\eqg}{\mbox{\thinspace\thinspace{\ms =}\thinspace\thinspace}}
\font\mss=msam8 at 9pt
\newcommand{\eqls}{\mbox{{\mss 5}}}

\font\msm=msbm8 at 12pt
\newcommand{\msb}[3]{\mbox{\msm \char'#1#2#3}}
\font\cm=cmssbx10 at 12pt
\font\cmb=cmssbx10 at 14pt
\begin{document}
\baselineskip=6mm
%
%
%
%
\begin{center}
{\cmb ON INVERSION OF \hspace{-1.3mm}{\large\boldmath{$H$}}-TRANSFORM
IN 
$\eufb114_{\nu ,r}$-SPACE}\\[7mm]
{\cm SERGEI A. SHLAPAKOV}\\[2mm]
Department of Mathematics\\
Vitebsk Pedagogical University, Vitebsk 210036, Belarus\\[4mm]
{\cm MEGUMI SAIGO}\\[2mm]
Department of Applied Mathematics\\
Fukuoka University, Fukuoka 814-0180, Japan\\[4mm]
{\cm ANATOLY A. KILBAS}\\[2mm]
Department of Mathematics and Mechanics\\
Belarusian State University, Minsk 220050, Belarus\\[10mm]
\end{center}
{\bf ABSTRACT.} \ The paper is devoted to study the inversion of the integral 
transform 
$$(\mbox{\boldmath$H$}f)(x)=\int^\infty_0H^{m,n}_{\thinspace p,q}
\left[xt\left|\begin{array}{c}(a_i,\alpha_i)_{1,p}\\[1mm](b_j,\beta_j)_{1,q}
\end{array}\right.\right]f(t)dt$$
involving the $H$-function as the kernel in the space $\euf114_{\nu ,r}$ of 
functions $f$ such that
$$\int^\infty_0\left|t^\nu f(t)\right|^r\frac{dt}t<\infty\quad(1<r<\infty,
\ \nu\in\msb122).$$
\\[2mm]
{\bf 1991 AMS SUBJECT CLASSIFICATION CODES:} \ 33C40, 44A20\\[3mm]
\begin{flushleft}
{\bf 1. \ INTRODUCTION}
\end{flushleft}
\setcounter{section}{1}
\par
This paper deals with the integral transforms of the form
$$(\mbox{\boldmath$H$}f)(x)=\int^\infty_0H^{m,n}_{\thinspace p,q}
\left[xt\left|\begin{array}{c}(a_i,\alpha_i)_{1,p}\\[2mm](b_j,\beta_j)_{1,q}
\end{array}\right.\right]f(t)dt,\eqno(1.1)$$
where $ H^{m,n}_{\thinspace p,q}\left[z\left|\begin{array}{c}
(a_i,\alpha_i)_{1,p}\\[2mm](b_j,\beta_j)_{1,q}\end{array}\right.\right]$ is 
the $H$-function, which is a function of general hypergeometric type being 
introduced by S. Pincherle in 1888 (see [2, \S1.19]). For integers $m,n,p,q$ 
such that $0\eql m\eql q,$ $0\eql n\eql p$, $a_i,b_j\in\msb103$ and 
$\alpha_i,\beta_j\in\msb122_+=[0,\infty)\ (1\eql i\eql p,1\eql j\eql q)$, 
it can be written by
\setcounter{equation}{1}
\begin{eqnarray}
H^{m,n}_{\thinspace p,q}\left[z\left|\begin{array}{c}(a_i,\alpha_i)_{1,p}
\\[2mm](b_j,\beta_j)_{1,q}\end{array}\right.\right]&\hspace{-2.5mm}=
&\hspace{-2.5mm}H^{m,n}_{\thinspace p,q}\left[z\left|\begin{array}{ll}
(a_1,\alpha_1),&\cdots,(a_p,\alpha_p)\\[2mm](b_1,\beta_1),&\cdots,
(b_q,\beta_q)\end{array}\right.\right]\nonumber\\[2mm]
&\hspace{-2.5mm}=&\hspace{-2.5mm}{\frac1{2\pi i}}\int_{L}
\eus110^{m,n}_{\thinspace p,q}\left[\left.\begin{array}{c}
(a_i,\alpha_i)_{1,p}\\[2mm](b_j,\beta_j)_{1,q}\end{array}\right|s\right]
z^{-s}ds,
\end{eqnarray}
where
$$\eus110^{m,n}_{\thinspace p,q}\left[\left.\begin{array}{c}
(a_i,\alpha_i)_{1,p}\\[2mm](b_j,\beta_j)_{1,q}\end{array}\right|s\right]=
{\frac{\displaystyle{\prod^m_{j=1}\Gamma(b_j+\beta_js)\prod^n_{i=1}
\Gamma(1-a_i-\alpha_is)}}{\displaystyle{\prod^p_{i=n+1}\Gamma(a_i+\alpha_is)
\prod^q_{j=m+1}\Gamma(1-b_j-\beta_js)}}},\eqno(1.3)$$
the contour $L$ is specially chosen and an empty product, if it occurs, is 
taken to be one. The theory of this function may be found in Braaksma [1], 
Srivastava {\it et al.} [13, Chapter 1], Mathai and Saxena [8, Chapter 2] 
and Prudnikov {\it et al.} [9, \S8.3]. We abbreviate the $H$-function (1.2) 
and the function (1.3) to $H(z)$ and $\eus110(s)$ when no confusion occurs. 
We note that the formal Mellin transform $\euf115$ of (1.1) gives the relation 
$$(\euf115\mbox{\boldmath$H$}f)(s)=\eus110(s)(\euf115f)(1-s).\eqno(1.4)$$
\par
Most of the known integral transforms can be put into the form (1.1), in 
particular, if $\alpha_{1}=\cdots =\alpha_{p}=\beta_{1}=\cdots =\beta_{q}=1$, 
(1.1) is the integral transform with Meijer's $G$-function in the kernel 
(Rooney [11], Samko {\it et al.} [12, \S36]). The integral transform (1.1) 
with the $H$-function kernel or the \mbox{\boldmath$H$}-transform was 
investigated by many authors (see Bibliography in Kilbas {\it et al.} [5-6]). 
In Kilbas {\it et al.} [5-7] we have studied it in the space $\euf114_{\nu ,r}
\ (1\eql r<\infty,\ \nu \in\msb122)$ consisted of Lebesgue measurable complex 
valued functions $f$ for which
$$\int^{\infty}_{0}\left|t^{\nu }f(t)\right|^{r}{\frac{dt}{t}}<\infty.
\eqno(1.5)$$
We have investigated the mapping properties such as the boundedness, the 
representation and the range of the \mbox{\boldmath$H$}-transform (1.1) on 
the space $\euf114_{\nu ,2}$ in Kilbas {\it et al.} [5] and on the space 
$\euf114_{\nu ,r}$ with any $1\eql r<\infty$ in Kilbas {\it et al.} [6-7], 
provided that $a^{\ast}\eqg 0,\ \delta =1$ and $\Delta =0$ or $\Delta \neq 0$, 
respectively. In Glaeske {\it et al.} [3] the results were extended to any 
$\delta>0$. Here
\setcounter{equation}{5}
\begin{eqnarray}
a^{\ast}&\hspace{-2.5mm}=&\hspace{-2.5mm}\sum^{n}_{i=1}\alpha_{i}-
\sum^{p}_{i=n+1}\alpha_{i}+\sum^{m}_{j=1}\beta_{j}-\sum^{q}_{j=m+1}\beta_{j};
\\[2mm]
\delta&\hspace{-2.5mm}=&\hspace{-2.5mm}\prod^{p}_{i=1}\alpha_{i}^{-\alpha_{i}}
\prod^{q}_{j=1}\beta_{j}^{\beta_{j}};\\[2mm]
\Delta&\hspace{-2.5mm}=&\hspace{-2.5mm}\sum^{q}_{j=1}\beta_{j}-\sum^{p}_{i=1}
\alpha_{i}.
\end{eqnarray}
In particular, we have proved that for certain ranges of parameters, the 
\mbox{\boldmath$H$}-transform (1.1) have the representations
\begin{eqnarray}
(\mbox{\boldmath$H$}f)(x)=hx^{1-(\lambda+1)/h}\displaystyle{\frac{d}{dx}}
x^{(\lambda+1)/h}\int^\infty_0H^{m,n+1}_{\thinspace p+1,q+1}\left[xt\left|
\begin{array}{l}(-\lambda,h),(a_i,\alpha_i)_{1,p}\\[2mm](b_j,\beta_j)_{1,q},
(-\lambda-1,h)\end{array}\right.\right]f(t)dt
\end{eqnarray}
or
\begin{eqnarray}
(\mbox{\boldmath$H$}f)(x)=-hx^{1-(\lambda+1)/h}\displaystyle{\frac{d}{dx}}
x^{(\lambda+1)/h}\!\int^\infty_0H^{m+1,n}_{\thinspace p+1,q+1}\left[xt
\left|\begin{array}{l}(a_i,\alpha_i)_{1,p},(-\lambda,h)\\[2mm](-\lambda-1,h),
(b_j,\beta_j)_{1,q}\end{array}\right.\right]f(t)dt,
\end{eqnarray}
owing to the value of ${\rm Re}(\lambda)$, where $\lambda\in\msb103$ and 
$h\in\msb122\setminus\{0\}$.
\par
In this paper we apply the results of Kilbas {\it et al.} [5-7] and Glaeske 
{\it et al.} [3] to find the inverse of the integral transforms (1.1) on the 
space $\euf114_{\nu ,r}$ with $1<r<\infty$ and $\nu\in\msb122$. Section 2 
contains preliminary information concerning the properties of the 
\mbox{\boldmath$H$}-transform (1.1) in the space $\euf114_{\nu ,r}$ and an 
asymptotic behavior of the $H$-function (1.2) at zero and infinity. In 
Sections 3 and 4 we prove that the inversion of the 
\mbox{\boldmath$H$}-transform have the respective form (1.9) or (1.10):
\begin{eqnarray}
&&f(x)=hx^{1-(\lambda+1)/h}\displaystyle{\frac{d}{dx}}x^{(\lambda+1)/h}
\nonumber\\[2mm]
&&\hspace{-18mm}\cdot\ \int^\infty_0H^{q-m,p-n+1}_{\thinspace p+1,q+1}
\left[xt\left|\begin{array}{l}(-\lambda,h),(1-a_i-\alpha_i,\alpha_i)_{n+1,p},
(1-a_i-\alpha_i,\alpha_i)_{1,n}\\[2mm](1-b_j-\beta_j,\beta_j)_{m+1,q},
(1-b_j-\beta_j,\beta_j)_{1,m},(-\lambda-1,h)\end{array}\right.\right]
(\mbox{\boldmath$H$}f)(t)dt
\end{eqnarray}
or
\begin{eqnarray}
&&f(x)=-hx^{1-(\lambda+1)/h}\displaystyle{\frac{d}{dx}}x^{(\lambda+1)/h}
\nonumber\\[2mm]
&&\hspace{-18mm}\cdot\ \int^\infty_0H^{q-m+1,p-n}_{\thinspace p+1,q+1}
\left[xt\left|\begin{array}{l}(1-a_i-\alpha_i,\alpha_i)_{n+1,p},
(1-a_i-\alpha_i,\alpha_i)_{1,n},(-\lambda,h)\\[2mm](-\lambda-1,h),
(1-b_j-\beta_j,\beta_j)_{m+1,q},(1-b_j-\beta_j,\beta_j)_{1,m}\end{array}
\right.\right](\mbox{\boldmath$H$}f)(t)dt,
\end{eqnarray}
provided that $a^{\ast}=0$. Section 3 is devoted to treat on the spaces 
$\euf114_{\nu ,2}$ and $\euf114_{\nu ,r}$ with $\Delta=0$, while Section 4 
on the space $\euf114_{\nu ,r}$ with $\Delta \neq 0$. 
\par
The obtained results are extensions of those by Rooney [11] from 
\mbox{\boldmath$G$}-transforms to \mbox{\boldmath$H$}-transforms.
\\
\begin{flushleft}
{\bf 2. \ PRELIMINARIES}
\end{flushleft}
\setcounter{section}{2}
\par
We give here some results from Kilbas {\it et al.} [5-6], Glaeske {\it et al.} 
[3] and from Kilbas and Saigo [4], Mathai and Saxena [8], Srivastava {\it et 
al.} [13] concerning the properties of \mbox{\boldmath$H$}-transforms (1.1) 
in $\euf114_{\nu ,r}$-spaces and the asymptotic behavior of the $H$-function 
at zero and infinity, respectively.
\par
For the $H$-function (1.2), let $a^{\ast}$ and $\Delta$ be defined by (1.6) 
and (1.8) and let
\setcounter{equation}{0}
\begin{eqnarray}
\alpha&\hspace{-2.5mm}=&\hspace{-2.5mm}\left\{\begin{array}{ll}\max
\left[-{\rm Re}\left(\displaystyle{\frac{b_1}{\beta_1}}\right),\cdots,
-{\rm Re}\left(\displaystyle{\frac{b_m}{\beta_m}}\right)\right]  & {\rm if}
\ m>0, \\\\ -\infty & {\rm if}\ m=0;\end{array}\right.\\[2mm]
\beta&\hspace{-2.5mm}=&\hspace{-2.5mm}\left\{\begin{array}{ll}\min
\left[{\rm Re}\left(\displaystyle{\frac{1-a_1}{\alpha_1}}\right),\cdots,
{\rm Re}\left(\displaystyle{\frac{1-a_n}{\alpha_n}}\right)\right]  & 
{\rm if}\ n>0, \\\\ \infty & {\rm if}\ n=0;\end{array}\right.\\[2mm]
a^{\ast}_1&\hspace{-2.5mm}=&\hspace{-2.5mm}\sum^m_{j=1}\beta_j-\sum^p_{i=n+1
}
\alpha_i;\quad a^{\ast}_2=\sum^n_{i=1}\alpha_i-\sum^q_{j=m+1}\beta_j;
\quad a^{\ast}_1+a^{\ast}_2=a^{\ast};\\[2mm]
\mu&\hspace{-2.5mm}=&\hspace{-2.5mm}\sum^{q}_{j=1}b_j-\sum^{p}_{i=1}a_i+
{\frac{p-q}{2}}.
\end{eqnarray}
\par
For the function $\eus110(s)$ given in (1.3), the exceptional set of 
$\eus110$ is meant the set of real numbers $\nu$ such that $\alpha <1-\nu 
<\beta$ and $\eus110(s)$ has a zero on the line ${\rm Re}(s)=1-\nu$ (see 
Rooney [11]). For two Banach space $X$ and $Y$ we denote by $[X,Y]$ the 
collection of bounded linear operators from $X$ to $Y$. 
\\\par
{\bf THEOREM 2.1.} \ [5, Theorem 3], [6, Theorem 3.3] \ {\sl Suppose that 
$\alpha<1-\nu<\beta$ and that either $a^{\ast}>0$ or $a^{\ast}=0,
\ \Delta(1-\nu)+{\rm Re}(\mu )\eql 0$. Then
\par
{\bf (a)} \ There is a one-to-one transform $\mbox{\boldmath$H$}\in
[\euf114_{\nu ,2},\euf114_{1-\nu ,2}]$ so that $(1.4)$ holds for $f\in
\euf114_{\nu ,2}$ and ${\rm Re}(s)=1-\nu$. If $a^{\ast}=0,\ \Delta 
(1-\nu)+{\rm Re}(\mu )=0$ and $\nu$ is not in the exceptional set of 
$\eus110,$ then the operator \mbox{\boldmath$H$} transforms $\euf114_{\nu ,2}$ 
onto $\euf114_{1-\nu ,2}$. 
\par
{\bf (b)} \ If $f\in \euf114_{\nu ,2}$ and ${\rm Re}(\lambda )>(1-\nu )h-1,$ 
$\mbox{\boldmath$H$}f$ is given by $(1.9).$ If $f\in \euf114_{\nu ,2}$ and 
${\rm Re}(\lambda )<(1-\nu )h-1,$ then $\mbox{\boldmath$H$}f$ is given by} 
(1.10). 
\\\par
{\bf THEOREM 2.2.} \ [6, Theorem 4.1], [3, Theorem 1] \ {\sl Let 
$a^{\ast}=\Delta =0,{\rm Re}(\mu)=0$ and $\alpha<1-\nu<\beta$.
\par
{\bf (a)} \ The transform $\mbox{\boldmath$H$}$ is defined on 
$\euf114_{\nu,2}$ and it can be extended to $\euf114_{\nu,r}$ as an element 
of $[\euf114_{\nu,r},\euf114_{1-\nu,r}]$ for $1<r<\infty$.
\par
{\bf (b)} \ If $1<r\eql 2,$ the transform $\mbox{\boldmath$H$}$ is one-to-one 
on $\euf114_{\nu ,r}$ and there holds the equality
$$(\euf115\mbox{\boldmath$H$}f)(s)=\eus110(s)(\euf115f)(1-s),\quad
{\rm Re}(s)=1-\nu.\eqno(2.5)$$
\par
{\bf (c)} \ If $f\in \euf114_{\nu ,r}\ (1<r<\infty),$ then 
$\mbox{\boldmath$H$}f$ is given by $(1.9)$ for ${\rm Re}(\lambda )>
(1-\nu )h-1,$ while $\mbox{\boldmath$H$}f$ is given by $(1.10)$ for} 
${\rm Re}(\lambda )<(1-\nu )h-1.$
\\\par
{\bf THEOREM 2.3.} \ [6, Theorem 5.1], [3, Theorem 3] \ {\sl Let $a^{\ast}=0,
\Delta>0,-\infty<\alpha <1-\nu <\beta,1<r<\infty$ and $\Delta (1-\nu)+
{\rm Re}(\mu)\eql 1/2-\gamma (r),$ where
$$\gamma (r)=\max\left[\frac1r,\frac1{r'}\right]\quad\mbox{with}\quad\frac1r+
\frac1{r'}=1.\eqno(2.6)$$
\par
{\bf (a)} \ The transform $\mbox{\boldmath$H$}$ is defined on 
$\euf114_{\nu ,2},$ and it can be extended to $\euf114_{\nu ,r}$ as an 
element of $[\euf114_{\nu ,r},\euf114_{1-\nu ,s}]$ for all $s$ with 
$r\eql s<\infty$ such that $s^{\prime}\eqg[1/2-\Delta(1-\nu)-{\rm Re}
(\mu)]^{-1}$ with $1/s+1/s^{\prime}=1$.
\par
{\bf (b)} \ If $1<r\eql 2,$ the transform $\mbox{\boldmath$H$}$ is one-to-one 
on $\euf114_{\nu ,r}$ and there holds the equality $(2.5).$
\par
{\bf (c)} \ If $f\in \euf114_{\nu ,r}$ and $g\in \euf114_{\nu ,s}$ with 
$1<r<\infty,1<s<\infty,1/r+1/s\eqg1$ and $\Delta(1-\nu)+{\rm Re}(\mu)\eql
1/2-\max [\gamma (r),\gamma (s)],$ then the relation
$$\int^{\infty}_{0}f(x)(\mbox{\boldmath$H$}g)(x)dx=\int^{\infty}_{0}g(x)
(\mbox{\boldmath$H$}f)(x)dx\eqno(2.7)$$
holds.}
\\\par
The following two assertions give the asymptotic behavior of the the 
$H$-function (1.2) at zero and infinity provided that the poles of Gamma 
functions in the numerator of $\eus110(s)$ do not coincide, i.e.
$$\beta_{j}(a_{i}-1-k)\neq \alpha_{i}(b_{j}+l)\quad
(i=1,\cdots,n;j=1,\cdots,m;k,l=0,1,2,\cdots).\eqno(2.8)$$ 
\\\par
{\bf THEOREM 2.4.} \ [8, \S1.1.6], [13, \S2.2] \ {\sl Let the condition 
$(2.8)$ be satisfied and poles of Gamma functions $\Gamma(b_j+\beta_js)
\ (j=1,\cdots ,m)$ be simple$,$ 
i.e.
$$\beta_{i}(b_{j}+k)\neq\beta_{j}(b_{i}+l)\quad
(i\neq j;i,j=1,\cdots ,m;k,l=0,1,2,\cdots).\eqno(2.9)$$  
If $\Delta\eqg 0,$ then
$$H^{m,n}_{p,q}(z)=O(z^{\rho})\quad(|z|\to 0)\quad\mbox{with}\quad\rho=
\min_{1\eqls j\eqls m}\left[{\frac{{\rm Re}(b_{j})}{\beta_j}}\right].
\eqno(2.10)$$
}
\\\par
{\bf THEOREM 2.5.} \ [4, Corollary 3] \ {\sl Let $a^{\ast},\Delta$ and $\mu$ 
be given by $(1.6), (1.8)$ and $(2.4),$ respectively. Let the conditions in 
$(2.8)$ be satisfied and poles of Gamma functions $\Gamma(1-a_i-\alpha_is)
\ (i=1,\cdots ,n)$ be simple$,$ 
i.e.
$$\alpha_{j}(1-a_{i}+k)\neq \alpha_{i}(1-a_{j}+l)\quad
(i\neq j;\ i,j=1,\cdots ,n;\ k,l=0,1,2 ,\cdots).\eqno(2.11)$$
If $a^{\ast}=0$ and $\Delta >0,$ then
$$H^{m,n}_{p,q}(z)=O(z^{\varrho })\quad(|z|\to\infty)\quad\mbox{with}
\quad\varrho =\max \left[\max_{1\eqls i\eqls n}\left[{\frac{{\rm Re}(a_{i})-1}
{\alpha_{i}}}\right],\ {\frac{{\rm Re}(\mu)+1/2}{\Delta }}\right].\eqno(2.12)$$
}
\\\par  
{\bf REMARK 2.1.} \ It was proved in Kilbas and Saigo [4, \S6] that if poles 
of Gamma functions $\Gamma(1-a_i-\alpha_is)\ (i=1,\cdots,n)$ are not simple 
(i.e. conditions in (2.11) are not satisfied), then the $H$-function (1.1) 
have power-logarithmic asymptotics at infinity. In this case the logaritmic 
multiplier $[\log (z)]^{N}$ with $N$ being the maximal number of orders of 
the poles may be added to the power multiplier $z^{\varrho}$ and hence the 
asymptotic estimate $O\left(z^{\varrho}\right)$ in (2.12) may be replaced by 
$O\left(z^{\varrho}[\log (z)]^{N}\right)$. The same result is valid in the 
case of the asymptotics of the $H$-function (1.1) at zero, and the estimate 
$O\left(z^{\rho}\right)$ in (2.10) may be replaced by $O\left(z^{\rho}
[\log (z)]^{M}\right)$, where $M$ is the maximal number of orders of the 
points at which the poles of $\Gamma(b_j+\beta_js)\ (j=1,\cdots,m)$ coincide.
\\\\
\begin{flushleft}
{\bf 3. INVERSION OF {\boldmath$H$}-TRANSFORM IN $\euf114_{\nu,2}$ AND 
$\euf114_{\nu,r}$ WHEN {\boldmath$\Delta$}$\thinspace={\bf0}$}
\end{flushleft}
\setcounter{section}{3}
\par
In this and next sections we investigate that $\mbox{\boldmath$H$}$-transform 
will have the inverse of the form (1.11) or (1.12). If 
$f\in \euf114_{\nu ,2}$, and $\mbox{\boldmath$H$}$ is defined on 
$\euf114_{\nu ,r}$, then according to Theorem 2.2, the equality (2.5) 
holds under the assumption there. This fact implies the relation 
$$(\euf115f)(s)={\frac{(\euf115\mbox{\boldmath$H$}f)(1-s)}{\eus110(1-s)}}
\eqno(3.1)$$
for ${\rm Re}(s)=\nu$. By (1.3) we have
$${\frac{1}{\eus110(1-s)}}=\eus110^{q-m,p-n}_{\thinspace p,q}
\left[\left.\begin{array}{c}(1-a_i-\alpha_i,\alpha_i)_{n+1,p},
(1-a_i-\alpha_i,\alpha_i)_{1,n}\\[2mm](1-b_j-\beta_j,\beta_j)_{m+1,q},
(1-b_j-\beta_j,\beta_j)_{1,m}\end{array}\right|s\right]\equiv\eus110_0(s),
\eqno(3.2)$$
and hence (3.1) takes the form
$$(\euf115f)(s)=(\euf115\mbox{\boldmath$H$}f)(1-s)\eus110_{0}(s)
\quad({\rm Re}(s)=\nu).\eqno(3.3)$$
We denote by $\alpha_{0},\beta_{0},{a}^{\ast}_{0},{a}^{\ast}_{01},
{a}^{\ast}_{02},\delta_{0},\Delta_{0}$ and $\mu_{0}$ for $\eus110_{0}$ 
instead of those for $\eus110$. Then we find
\setcounter{equation}{3}
\begin{eqnarray}
&&\alpha_{0}=\left\{\begin{array}{ll}\max\left[\displaystyle{\frac{{\rm Re}
(b_{m+1})-1}{\beta_{m+1}}}+1,\cdots,\displaystyle{\frac{{\rm Re}(b_{q})-1}
{\beta_{q}}}+1\right]  & {\rm if}\ q>m, \\\\ -\infty & {\rm if}\ q=m;
\end{array}\right.\\[2mm]
&&\beta_{0}=\left\{\begin{array}{ll}\min\left[\displaystyle{\frac{{\rm Re}
\left(a_{n+1}\right)}{\alpha_{n+1}}+1},\cdots,\displaystyle{\frac{{\rm Re}
\left(a_{p}\right)}{\alpha_{p}}}+1\right]  & {\rm if}\ p>n, \\\\ \infty & 
{\rm if}\ p=n;\end{array}\right.\\[2mm]
&&{a}^{\ast}_{0}=-{a}^{\ast};\ \ {a}^{\ast}_{01}=-{a}^{\ast}_{2};
\ \ {a}^{\ast}_{02}=-{a}^{\ast}_{1};\ \ \delta_{0}=\delta;
\ \ \Delta_{0}=\Delta;\ \ \mu_{0}=-\mu -\Delta.
\end{eqnarray}
We also note that if $\alpha_{0}<\nu <\beta_{0}$, $\nu$ is not in the 
exceptional set of $\eus110_{0}$.
\par
First we consider the case $r=2$.
\\\par
{\bf THEOREM 3.1.} \ {\sl Let $\alpha<1-\nu<\beta,\alpha_{0}<
\nu <\beta_{0},a^{\ast}=0$ and $\Delta(1-\nu)+{\rm Re}(\mu)=0$. If $f\in
\euf114_{\nu ,2},$ the relation $(1.11)$ holds for ${\rm Re}(\lambda)>\nu h-1$ 
and the relation $(1.12)$ holds for} ${\rm Re}(\lambda)<\nu h-1$.
\par
{\bf PROOF.} \ We apply Theorem 2.1 with $\eus110$ being replaced by 
$\eus110_{0}$ and $\nu$ by $1-\nu$. By the assumption and (3.6) we have
\begin{eqnarray}
&&{a}^{\ast}_{0}=-{a}^{\ast}=0,\\[2mm]
&&\Delta_{0}[1-(1-\nu)]+{\rm Re}(\mu_{0})=\Delta \nu -{\rm Re}(\mu)-\Delta=
-[\Delta(1-\nu)+{\rm Re}(\mu)]=0
\end{eqnarray}
and $\alpha_{0}<1-(1-\nu)<\beta_{0}$, and thus Theorem 2.1(a) applies. Then 
there is a one-to-one transform $\mbox{\boldmath$H$}_{0}\in[\euf114_{1-\nu,2},
\euf114_{\nu,2}]$ so that the relation
$$(\euf115\mbox{\boldmath$H$}_{0}f)(s)=\eus110_{0}(s)(\euf115f)(1-s)
\eqno(3.9)$$
holds for $f\in \euf114_{1-\nu,2}$ and ${\rm Re}(s)=\nu$. Further if $f\in
\euf114_{\nu,2}$, $\mbox{\boldmath$H$}f\in \euf114_{1-\nu,2}$ and it follows 
from (3.9), (1.4) and (3.2) that
$$(\euf115\mbox{\boldmath$H$}_{0}\mbox{\boldmath$H$}f)(s)=
\eus110_{0}(s)(\euf115\mbox{\boldmath$H$}f)(1-s)=\eus110_{0}(s)
\eus110(1-s)(\euf115f)(s)=(\euf115f)(s),$$
if ${\rm Re}(s)=\nu$. Hence $\euf115\mbox{\boldmath$H$}_{0}
\mbox{\boldmath$H$}f=\euf115f$ and 
$$\mbox{\boldmath$H$}_{0}\mbox{\boldmath$H$}f=f\quad\mbox{for}\quad f\in
\euf114_{\nu,2}.\eqno(3.10)$$
\par
Applying Theorem 2.1(b) with $\eus110$ being replaced by $\eus110_{0}$ and 
$\nu$ by $1-\nu$, we obtain for $f\in \euf114_{1-\nu,2}$ that
\setcounter{equation}{10}
\begin{eqnarray}
&&\hspace{-15mm}(\mbox{\boldmath$H$}_{0}f)(x)=hx^{1-(\lambda+1)/h}
\displaystyle{\frac{d}{dx}}x^{(\lambda+1)/h}\nonumber\\[2mm]
&&\hspace{-12mm}\cdot\ \int^\infty_0H^{q-m,p-n+1}_{\thinspace p+1,q+1}
\left[xt\left|\begin{array}{l}(-\lambda,h),(1-a_i-\alpha_i,\alpha_i)_{n+1,p},
(1-a_i-\alpha_i,\alpha_i)_{1,n}\\[2mm](1-b_j-\beta_j,\beta_j)_{m+1,q},
(1-b_j-\beta_j,\beta_j)_{1,m},(-\lambda-1,h)\end{array}\right.\right]f(t)dt,
\end{eqnarray}
if ${\rm Re}(\lambda)>[1-(1-\nu )]h-1$ and 
\begin{eqnarray}
&&\hspace{-15mm}(\mbox{\boldmath$H$}_{0}f)(x)=-hx^{1-(\lambda+1)/h}
\displaystyle{\frac{d}{dx}}x^{(\lambda+1)/h}\nonumber\\[2mm]
&&\hspace{-12mm}\cdot\ \int^\infty_0H^{q-m+1,p-n}_{\thinspace p+1,q+1}
\left[xt\left|\begin{array}{l}(1-a_i-\alpha_i,\alpha_i)_{n+1,p},
(1-a_i-\alpha_i,\alpha_i)_{1,n},(-\lambda,h)\\[2mm](-\lambda-1,h),
(1-b_j-\beta_j,\beta_j)_{m+1,q},(1-b_j-\beta_j,\beta_j)_{1,m}\end{array}
\right.\right]f(t)dt,
\end{eqnarray}
if ${\rm Re}(\lambda)<[1-(1-\nu )]h-1$. Replacing $f$ by 
$\mbox{\boldmath$H$}f$ and using (3.10) we have the relations (1.11) and 
(1.12) for $f\in \euf114_{\nu,2}$, if ${\rm Re}(\lambda)>\nu h-1$ and 
${\rm Re}(\lambda)<\nu h-1$, respectively, which completes the proof of 
theorem.
\\\par
Next results is the extension of Theorem 3.1 to $\euf114_{\nu,r}$-spaces for 
any $1<r<\infty,$ provided that $\Delta =0$ and ${\rm Re}(\mu)=0$. 
\\\par
{\bf THEOREM 3.2.} \ {\sl Let $\alpha <1-\nu <\beta,\alpha_{0}<\nu
<\beta_{0},a^{\ast}=0,\Delta =0$ and ${\rm Re}(\mu)=0$. If $f\in
\euf114_{\nu ,r}\ (1<r<\infty),$ the relation $(1.11)$ holds for 
${\rm Re}(\lambda)>\nu h-1$ and the relation $(1.12)$ holds for} 
${\rm Re}(\lambda)<\nu h-1$.
\par
{\bf PROOF.} \ We apply Theorem 2.2 with $\eus110$ being replaced by 
$\eus110_0$ and $\nu$ by $\nu-1$. By the assumption and (3.6), we have 
$a^*_{0}=\Delta_{0}=0,{\rm Re}(\mu_{0})=0$ and $\alpha_{0}<1-(1-\nu)<
\beta_{0}$, and thus Theorem 2.2(a) can be applied. In accordance with this 
theorem, $\mbox{\boldmath$H$}_{0}$ can be extended to $\euf114_{1-\nu,r}$ as 
an element of $\mbox{\boldmath$H$}_{0}\in [\euf114_{1-\nu,r},
\euf114_{\nu,r}]$. By virtue of (3.10) $\mbox{\boldmath$H$}_{0}
\mbox{\boldmath$H$}$ is identical operator in $\euf114_{\nu,2}$. By Rooney 
[11, Lemma 2.2] $\euf114_{\nu,2}$ is dense in $\euf114_{\nu,r}$ and since 
$\mbox{\boldmath$H$}\in [\euf114_{\nu,r},\euf114_{1-\nu,r}]$ and 
$\mbox{\boldmath$H$}_{0}\in [\euf114_{1-\nu,r},\euf114_{\nu,r}]$, the 
operator $\mbox{\boldmath$H$}_{0}\mbox{\boldmath$H$}$ is identical in 
$\euf114_{\nu,r}$ and hence 
$$\mbox{\boldmath$H$}_{0}\mbox{\boldmath$H$}f=f\quad\mbox{for}\quad f\in
\euf114_{\nu,r}.\eqno(3.13)$$
\par
Applying Theorem 2.2(c) with $\eus110$ being replaced by $\eus110_{0}$ and 
$\nu$ by $1-\nu$, we obtain that the relations (3.11) and (3.12) hold for 
$f\in \euf114_{1-\nu,r}$, when ${\rm Re}(\lambda)>[1-(1-\nu )]h-1$ and 
${\rm Re}(\lambda)<[1-(1-\nu )]h-1$, respectively. Replacing $f$ by 
$\mbox{\boldmath$H$}f$ and using (3.13), we arrive at (1.11) and (1.12) for 
$f\in \euf114_{1-\nu,r}$, if ${\rm Re}(\lambda)>\nu h-1$ and 
${\rm Re}(\lambda)<\nu h-1$, respectively, which completes the proof of 
theorem.
\\\par
{\bf REMARK 3.1.} \ If $\alpha_{1}=\cdots =\alpha_{p}=\beta_{1}=\cdots =
\beta_{q}=1$ which means that the $H$-function (1.2) is Meijer's $G$-function, 
then $\Delta=q-p$ and Theorems 8.1 and 8.2 in Rooney [11] follow from 
Theorems 3.1 and 3.2.
\\
\begin{flushleft}
{\bf 4. INVERSION OF {\boldmath$H$}-TRANSFORM IN $\euf114_{\nu ,r}$ WHEN 
{\boldmath$\Delta$}$\thinspace\ne{\bf0}$}
\end{flushleft}
\setcounter{section}{4}
\par
We now investigate under what condition the $\mbox{\boldmath$H$}$-transform 
with $\Delta \neq 0$ will have the inverse of the form (1.11) or (1.12). 
First, we consider the case $\Delta >0$. To obtain the inversion of the 
$\mbox{\boldmath$H$}$-transform on $\euf114_{\nu ,r}$ we use the relation 
(2.7).
\\\par
{\bf THEOREM 4.1.} \ {\sl Let $1<r<\infty,-\infty<\alpha<1-\nu <\beta,
\alpha_0<\nu<\min\{\beta_0,[{\rm Re}(\mu +1/2)/\Delta]+1\},a^{\ast}=0,
\Delta >0$ and $\Delta (1-\nu)+{\rm Re}(\mu)\eql 1/2-\gamma(r),$ where 
$\gamma(r)$ is given in $(2.6)$. If $f\in \euf114_{\nu ,r},$ then the 
relations $(1.11)$ and $(1.12)$ hold for ${\rm Re}(\lambda)>\nu h-1$ and 
for ${\rm Re}(\lambda)<\nu h-1,$ respectively.}
\par
{\bf PROOF.} \ According to Theorem 2.3(a), the 
$\mbox{\boldmath$H$}$-transform is defined on $\euf114_{\nu ,r}$. First we 
consider the case ${\rm Re}(\lambda)>\nu h-1$. Let $H_{1}(t)$ be the function
$$H_{1}(t)=H^{q-m,p-n+1}_{\thinspace p+1,q+1}\left[t\left|\begin{array}{l}
(-\lambda,h),(1-a_i-\alpha_i,\alpha_i)_{n+1,p},(1-a_i-\alpha_i,\alpha_i)_{1,n}
\\[2mm](1-b_j-\beta_j,\beta_j)_{m+1,q},(1-b_j-\beta_j,\beta_j)_{1,m},
(-\lambda-1,h)\end{array}\right.\right].\eqno(4.1)$$
If we denote by $\widetilde{{a}}^{\ast},\widetilde{\delta},\widetilde{\Delta}$ 
and $\widetilde{\mu}$ for $H_{1}$ instead of those for $H,$ then 
$$\widetilde{a}^{\ast}=-a^{\ast}=0;\quad\widetilde{\delta}=\delta;
\quad\widetilde{\Delta}=\Delta >0;\quad\widetilde{\mu}=-\mu -\Delta -1.
\eqno(4.2)$$
\par
We prove that $H_{1}\in \euf114_{\nu ,s}$ for any $s\ (1\eql s<\infty)$. For 
this, we first apply Theorems 2.4 and 2.5 and Remark 2.1 to $H_{1}(t)$ to 
find its asymptotic behavior at zero and infinity. According to (3.4), (3.5) 
and the assumptions, we find
\begin{eqnarray*}
{\frac{{\rm Re}(b_{j})-1}{\beta_{j}}}+1&\hspace{-2.5mm}\eql&\hspace{-2.5mm}
\alpha_{0}<\beta_{0}\eql{\frac{{\rm Re}(a_{i})}{\alpha_{i}}}+1\quad(j=m+1,
\cdots ,q;\ i=n+1,\cdots,p);\\[2mm]
{\frac{{\rm Re}(b_{j})-1}{\beta_{j}}}+1&\hspace{-2.5mm}\eql&\hspace{-2.5mm}
\alpha_{0}<\nu<{\frac{{\rm Re}(\lambda )+1}{h}}\quad(j=m+1,\cdots,q).
\end{eqnarray*}
Then it follows from here that the poles 
$$a_{ik}={\frac{a_{i}+k}{\alpha_{i}}}+1\quad(i=n+1,\cdots,p;\ k=0,1,2,\cdots),
\quad\lambda_{n}={\frac{\lambda +1+n}{h}}\quad(n=0,1,2,\cdots )$$ 
of Gamma functions $\Gamma(a_i+\alpha_i-\alpha_{i}s)\ (i=n+1,\cdots,p)$ and 
$\Gamma(1+\lambda-hs)$, and the poles 
$$b_{jl}={\frac{b_{j}-1-l}{\beta_{j}}}+1\quad(j=m+1,\cdots,q;\ l=0,1,2,
\cdots)$$
of Gamma functions $\Gamma(1-b_j-\beta_{j}+\beta_js)\ (j=m+1,\cdots ,q)$ do 
not coincide. Hence by Theorem 2.4, (4.1) and Remark 2.1, we have 
$$H_{1}(t)=O(t^{\rho_{1}})\quad(|t|\to 0)\quad\mbox{with}\quad\rho_{1} =
\min_{m+1\eqls j\eqls q}\left[{\frac{1-{\rm Re}(b_{j})}{\beta_{j}}}\right]-1=
-\alpha_{0}$$
for $\alpha_{0}$ being given in (3.4), or
$$H_{1}(t)=O(t^{-\alpha_{0}})\quad(t\to 0)\eqno(4.3)$$
with an additional logarithmic multilplier $[\log t]^N$ possibly, if Gamma 
functions $\Gamma(1-b_j-\beta_j+\beta_js)\ (j=m+1,\cdots ,q)$ have general 
poles of order $N\eqg 2$ at some points.
\par
Further by Theorem 2.5, (4.1) and Remark 2.1,
$$H_{1}(t)=O(t^{\varrho_{1}})\quad(t\to\infty)\quad\mbox{with}
\quad\varrho_{1} =\max \left[\beta_{0},{\frac{-{\rm Re}(\mu)-1/2}
{\Delta}}-1,{\frac{-{\rm Re}(\lambda)-1}{h}}\right]$$
for $\beta_{0}$ being given by (3.5), or
$$H_{1}(t)=O(t^{-\gamma_{0}})\quad(|t|\to\infty)\quad\mbox{with}
\quad\gamma_{0}=\min \left[\beta_{0},{\frac{{\rm Re}(\mu)+1/2}{\Delta}}+1,
{\frac{{\rm Re}(\lambda)+1}{h}}\right]\eqno(4.4)$$
and with an additional logarithmic multilplier $[\log (t)]^{M}$ possibly, 
if Gamma functions $\Gamma (1+\lambda -hs),\Gamma(a_i+\alpha_i-\alpha_is)
\ (i=n+1,\cdots ,p)$ have general poles of order $M\eqg 2$ at some points.
\par
Let Gamma functions $\Gamma(1-b_j-\beta+\beta_js)\ (j=m+1,\cdots ,q)$ and 
$\Gamma (1+\lambda -hs),\Gamma(a_i+\alpha_i-\alpha_is)\ (i=n+1,\cdots ,p)$ 
have simple poles. Then from (4.3) and (4.4) we see that for $1\eql s<\infty$, 
$H_{1}(t)\in \euf114_{\nu ,s}$ if and only if, for some $R_{1}$ and 
$R_{2},\ 0<R_{1}<R_{2}<\infty$, the integrals
$$\int^{R_{1}}_{0}t^{s(\nu-\alpha_{0})-1}dt,\quad\int^{\infty}_{R_{2}}
t^{s(\nu-\gamma_{0})-1}dt\eqno(4.5)$$
are convergent. Since by the assumption $\nu >\alpha_{0}$, the first integral 
in (4.5) converges. In view of our assumtions 
$$\nu <\beta_{0},\quad\nu <{\frac{{\rm Re}(\mu)+1/2}{\Delta}}+1,\quad\nu <
{\frac{{\rm Re}(\lambda)+1}{h}}$$ 
we find $\nu -\gamma_{0}<0$ and the second integral in (4.5) converges, too.
\par
If Gamma functions $\Gamma(1-b_j-\beta_{j}+\beta_js)\ (j=m+1,\cdots ,q)$ or 
$\Gamma (1+\lambda -hs),\Gamma(a_i+\alpha_{i}-\alpha_is)\ (i=n+1,\cdots ,p)$ 
have general poles, then the logarithmic multipliers $[\log (t)]^{N}
\ (N=1,2,\cdots)$ may be added in the integrals in (4.5), but they do not 
influence on the convergence of them. Hence, under the assumptions we have
$$H_{1}(t)\in \euf114_{\nu ,s}\quad(1\eql s<\infty).\eqno(4.6)$$
\par
Let $a$ be a positive number and $\Pi_{a}$ denote the operator
$$(\Pi_{a}f)(x)=f(ax)\quad(x>0)\eqno(4.7)$$
for a function $f$ defined almost everywhere on $(0,\infty)$. It is known in 
Rooney [11, p.268] that $\Pi_{a}$ is a bounded isomorphism of 
$\euf114_{\nu ,r}$ onto $\euf114_{a \nu ,r}$, and if $f\in \euf114_{\nu,r}
\ (1\eql r\eql 2),$ there holds the relation for the Mellin transform 
$\euf115$
$$(\euf115\Pi_{a}f)(s)=a^{-s}(\euf115f)\left(\frac sa\right)\quad({\rm Re}(s)=
\nu).\eqno(4.8)$$ 
\par
By virtue of Theorem 2.3(c) and (4.6), if $f\in \euf114_{\nu ,r}$ and 
$H_{1}\in \euf114_{\nu ,r^{\prime}}$ (and hence $\Pi_{x}H_{1}\in
\euf114_{\nu ,r^{\prime}}$), then 
$$\int^{\infty}_{0}H_{1}(xt)(\mbox{\boldmath$H$}f)(t)dt=\int^{\infty}_{0}
(\Pi_{x}H_{1})(t)(\mbox{\boldmath$H$}f)(t)dt=\int^{\infty}_{0}
(\mbox{\boldmath$H$}\Pi_{x}H_{1})(t)f(t)dt.\eqno(4.9)$$
From the assumption $\Delta (1-\nu)+{\rm Re}(\mu)\eql 1/2-\gamma (r)\eql 0$, 
Theorem 2.3(b) and (4.8) imply that
$$(\euf115\mbox{\boldmath$H$}\Pi_{x}H_{1})(s)=\eus110(s)(\euf115\Pi_{x}H_{1})
(1-s)=x^{-(1-s)}\eus110(s)(\euf115H_{1})(1-s)\eqno(4.10)$$
for ${\rm Re}(s)=1-\nu$. Now from (4.6), $H_{1}(t)\in \euf114_{\nu ,1}$. Then 
by the definitions of the $H$-function (1.2), (1.3) and the direct and inverse 
Mellin transforms (see, for example, Samko {\it et al.} [12, (1.112), 
(1.113)]), we have
\begin{eqnarray*}
(\euf115H_{1})(s)&\hspace{-2.5mm}=&\hspace{-2.5mm}
\eus110^{q-m,p-n+1}_{\thinspace p+1,q+1}\left[\left.\begin{array}{l}
(-\lambda,h),(1-a_i-\alpha_i,\alpha_i)_{n+1,p},(1-a_i-\alpha_i,\alpha_i)_{1,n}
\\[2mm](1-b_j-\beta_j,\beta_j)_{m+1,q},(1-b_j-\beta_j,\beta_j)_{1,m},
(-\lambda-1,h)\end{array}\right|s\right]\\[2mm]
&\hspace{-2.5mm}=&\hspace{-2.5mm}\eus110^{q-m,p-n}_{\thinspace p,q}
\left[\left.\begin{array}{l}(1-a_i-\alpha_i,\alpha_i)_{n+1,p},
(1-a_i-\alpha_i,\alpha_i)_{1,n}\\[2mm](1-b_j-\beta_j,\beta_j)_{m+1,q},
(1-b_j-\beta_j,\beta_j)_{1,m}\end{array}\right|s\right]
{\frac{\Gamma (1+\lambda -hs)}{\Gamma (2+\lambda -hs)}}\\[2mm]
&\hspace{-2.5mm}=&\hspace{-2.5mm}{\frac{\eus110_0(s)}{1+\lambda -hs}}
\end{eqnarray*}
for ${\rm Re}(s)=\nu$, where $\eus110_0$ is given by (3.2). It follows from 
here that for ${\rm Re}(s)=1-\nu$,
$$(\euf115H_{1})(1-s)={\frac{\eus110_0(1-s)}{1+\lambda -h(1-s)}}=
{\frac{1}{\eus110(s)[1+\lambda -h(1-s)]}}.$$
Substituting this into (4.10) we obtain
$$(\euf115\mbox{\boldmath$H$}\Pi_{x}H_{1})(s)={\frac{x^{-(1-s)}}
{1+\lambda -h(1-s)}}\quad({\rm Re}(s)=1-\nu).\eqno(4.11)$$
\par
For $x>0$ let us denote by $g_{x}(t)$ a function 
$$g_{x}(t)=\left\{\begin{array}{ll}\displaystyle{\frac 1h}
\ t^{(\lambda +1)/h-1}  & {\rm if}\ 0<t<x, \\\\0 & {\rm if}\ t>x, 
\end{array}\right. \eqno(4.12)$$
then 
$$(\euf115g_{x})(s)={\frac{x^{s+(\lambda +1)/h-1}}{1+\lambda -h(1-s)}},$$
and (4.11) takes the form
$$(\euf115\mbox{\boldmath$H$}\Pi_{x}H_{1})(s)=(\euf115[x^{-(\lambda +1)/h}
g_{x})])(s),$$
which implies
$$(\mbox{\boldmath$H$}\Pi_{x}H_{1})(t)=x^{-(\lambda +1)/h}g_{x}(t).
\eqno(4.13)$$
\par
Substituting (4.13) into (4.9), we have  
$$\int^{\infty}_{0}H_{1}(xt)(\mbox{\boldmath$H$}f)(t)dt=x^{-(\lambda +1)/h}
\int^{\infty}_{0}g_{x}(t)f(t)dt$$
or, in accordance with (4.12),
$$\int^{x}_{0}t^{(\lambda +1)/h-1}f(t)dt=hx^{(\lambda +1)/h}\int^{\infty}_{0}
H_{1}(xt)(\mbox{\boldmath$H$}f)(t)dt.$$
Differentiating this relation we obtain
$$f(x)=hx^{1-(\lambda +1)/h}{\frac{d}{dx}}x^{(\lambda +1)/h}\int^{\infty}_{0}
H_{1}(xt)(\mbox{\boldmath$H$}f)(t)dt$$
which shows (1.11).
\par
If ${\rm Re}(\lambda)<\nu h-1$, the relation (1.12) is proved similarly to 
(1.11), by taking the function
$$H_{2}(t)=H^{q-m+1,p-n}_{\thinspace p+1,q+1}\left[t\left|\begin{array}{l}
(1-a_i-\alpha_i,\alpha_i)_{n+1,p},(1-a_i-\alpha_i,\alpha_i)_{1,n},(-\lambda,h)
\\[2mm](-\lambda-1,h),(1-b_j-\beta_j,\beta_j)_{m+1,q},
(1-b_j-\beta_j,\beta_j)_{1,m}\end{array}\right.\right]\eqno(4.14)$$
instead of the function $H_{1}(t)$ in (4.1). This completes the proof of the 
theorem.
\\\par
In the case $\Delta <0$ the following statement gives the inversion of 
$\mbox{\boldmath$H$}$-transform on $\euf114_{\nu ,r}$.
\\
\par
{\bf THEOREM 4.2.} \ {\sl Let $1<r<\infty,\alpha <1-\nu <\beta<\infty,
\max[\alpha_{0},\{{\rm Re}(\mu +1/2)/\Delta\}+1]<\nu<\beta_{0},a^{\ast}=0,
\Delta <0$ and $\Delta (1-\nu)+{\rm Re}(\mu)\eql 1/2-\gamma (r),$ where 
$\gamma (r)$ is given by $(2.6).$ If $f\in \euf114_{\nu ,r},$ then the 
relations $(1.11)$ and $(1.12)$ holds for ${\rm Re}(\lambda)>\nu h-1$ and for 
${\rm Re}(\lambda)<\nu h-1,$ respectively.}
\\
\par
This theorem is proved similarly to Theorem 4.1, if we apply Theorem 5.2 
from Kilbas {\it et al.} [6] instead of Theorem 2.3 and take into account 
the asymptotics of the $H$-function at zero and infinity (see Srivastava 
{\it et al.} [13, \S2.2] and Kilbas and Saigo [4, Corollary 4]).
\\
\par
{\bf REMARK 4.1.} \ If $\alpha_{1}=\cdots =\alpha_{p}=\beta_{1}=\cdots =
\beta_{q}=1$, then Theorems 8.3 and 8.4 in Rooney [11] follow from 
Theorems 4.1 and 4.2.
\\
\begin{center}{\bf REFERENCES}\end{center}
\baselineskip=4.47mm
\begin{enumerate}
\item[{[1]\ }]BRAAKSMA, B.L.G. Asymptotic expansions and analytic continuation 
for a class of Barnes integrals, {\it Compos. Math.} {\bf 15}(1964), 239-341.
\item[{[2]\ }]ERD\'ELYI, A., MAGNUS, W., OBERHETTINGER, F. and TRICOMI, 
F.G. {\it Higher Transcendental Functions} Vol.\thinspace1, McGraw-Hill, New 
York-Toronto-London, 1953.
\item[{[3]\ }]GLAESKE, H.-J., KILBAS, A.A., SAIGO, M. and SHLAPAKOV S.A. 
$L_{\nu,r}$-theory of integral transforms with $H$-function as kernels 
(Russian), {\it Dokl. Akad. Nauk Belarusi} {\bf 41}(1997), 10-15.
\item[{[4]\ }]KILBAS, A.A. and SAIGO, M. On asymptotics of Fox's $H$-function 
at zero and infinity, {\it First International Workshop}, {\it Transform 
Methods and Special Functions}, (Sofia, Bulgaria, 1994), 99-122, Science 
Culture Technology Publ., Singapore, 1995.
\item[{[5]\ }]KILBAS, A.A., SAIGO, M. and SHLAPAKOV, S.A. Integral transforms 
with Fox's $H$-function in spaces of summable functions, {\it Integral Transf. 
Specc. Func.}, {\bf 1}(1993), 87-103.
\item[{[6]\ }]KILBAS, A.A., SAIGO, M. and SHLAPAKOV, S.A. Integral transforms 
with Fox's $H$-function in $\euf114_{\nu ,r}$-spaces, {\it Fukuoka Univ. Sci. 
Rep.} {\bf 23}(1993), 9-31.
\item[{[7]\ }]KILBAS, A.A., SAIGO, M. and SHLAPAKOV, S.A. Integral transforms 
with Fox's $H$-function in $\euf114_{\nu ,r}$-spaces. II, {\it Fukuoka Univ. 
Sci. Rep.} {\bf 24}(1994), 13-38.
\item[{[8]\ }]MATHAI, A.M. and SAXENA, R.K. {\it The $H$-Function with 
Applications in Statistics and other Disciplines}, Wiley Eastern, New Delhi, 
1978.
\item[{[9]\ }]PRUDNIKOV, A.P., BRYCHKOV, Yu.A. and MARICHEV, O.I. 
{\it Integrals and Series, Vol.3$:$ More Special Functions}, Gordon and 
Breach, New York {\it et alibi}, 1990.  
\item[{[10]\ }]ROONEY, P.G. A technique for studying the boundedness and 
extendability of certain types of operators. {\it Canad. J. Math.} 
{\bf 25}(1973), 1090-1102. 
\item[{[11]\ }]ROONEY, P.G. On integral transformations with $G$-function 
kernels, {\it Proc. Royal Soc. Edinburgh} {\bf A93}(1983), 265-297.
\item[{[12]\ }]SAMKO, S.G., KILBAS, A.A. and MARICHEV, O.I. {\it Fractional 
Integrals and Derivatives. Theory and Applications}, Gordon and Breach, 
Yverdon (Switzerlan) {\it et alibi}, 1993.
\item[{[13]\ }]SRIVASTAVA, H.M., GUPTA, K.C. and GOYAL, S.P. 
{\it The $H$-Functions of One and Two Variables with Applications}, South 
Asian Publishers, New Delhi-Madras, 1982.
\end{enumerate}
\end{document}